\documentclass[11pt,twoside]{article}
\usepackage{epsfig,amsfonts,color}
\usepackage{amsmath,bm}
\usepackage{mathrsfs}
\allowdisplaybreaks
\usepackage{multirow}
\usepackage{cite}
\usepackage{amssymb, palatino, geometry,url}
\usepackage[colorlinks=true,linkcolor=blue,citecolor=blue,urlcolor=blue]{hyperref}
\usepackage{lineno}
\usepackage{enumerate}
\usepackage{float}
\usepackage{graphicx}
\usepackage{array}
\usepackage{placeins}
\usepackage{caption}
\usepackage{subcaption}
\usepackage{pdflscape}
\usepackage{algorithm}
\usepackage{algpseudocode}
\usepackage{fancyhdr,listings}
\usepackage{enumitem}

\newcommand{\mynote}[1]{{\color{black}{#1}}}

\newcommand{\be}{\begin{equation}}
\newcommand{\ee}{\end{equation}}
\newcommand{\bea}{\begin{eqnarray}}
\newcommand{\eea}{\end{eqnarray}}

\newcommand{\bvec}{\left(\begin{array}{c}}
\newcommand{\evec}{\end{array}\right)}
\newcommand{\bsub}{\begin{subequations}}
\newcommand{\esub}{\end{subequations}}

\newcolumntype{L}[1]{>{\raggedright\let\newline\\\arraybackslash\hspace{0pt}}m{#1}}
\newcolumntype{C}[1]{>{\centering\let\newline\\\arraybackslash\hspace{0pt}}m{#1}}
\newcolumntype{R}[1]{>{\raggedleft\let\newline\\\arraybackslash\hspace{0pt}}m{#1}}

\lstdefinelanguage{julia}
{
keywordsprefix=\@,
morekeywords={
exit,whos,edit,load,is,isa,isequal,typeof,tuple,ntuple,uid,hash,finalizer,convert,promote,
subtype,typemin,typemax,realmin,realmax,sizeof,eps,promote_type,method_exists,applicable,
invoke,dlopen,dlsym,system,error,throw,assert,new,Inf,Nan,pi,im,begin,while,for,in,return,
break,continue,macro,quote,let,if,elseif,else,try,catch,end,bitstype,ccall,do,using,module,
import,export,importall,baremodule,immutable,local,global,const,Bool,Int,Int8,Int16,Int32,
Int64,Uint,Uint8,Uint16,Uint32,Uint64,Float32,Float64,Complex64,Complex128,Any,Nothing,None,
function,type,typealias,abstract
},
sensitive=true,
morecomment=[l]{\#},
morestring=[b]',
morestring=[b]"
}

\begin{document}


\title{Multiscale Model Predictive Control of \\ Battery Systems for Frequency Regulation Markets \\
 using Physics-Based Models}

\author{Yankai~Cao${}^{*}$, Seong~Beom~Lee${}^\ddag$, Venkat~R.~Subramanian${}^\dag$, Victor M. Zavala${}^{\P}$\\
{\small ${}^{*}$Department of Chemical and Biological Engineering}\\
{\small \;University of British Columbia, 5000--2332 Main Mall, Vancouver, BC, Canada}\\
{\small ${}^{\P}$Department of Chemical and Biological Engineering}\\
{\small \;University of Wisconsin-Madison, 1415 Engineering Dr, Madison, WI 53706, USA}\\
{\small ${}^\ddag$Department of Chemical Engineering}\\
{\small \;University of Washington, 3781 Okanogan Ln, Seattle, WA 98195, USA}\\
{\small ${}^\dag$Department of Mechanical Engineering}\\
{\small \;University of Texas at Austin, 04 E. Dean Keeton Street, SAustin, TX 78712, USA}
}
 \date{}
\maketitle

\begin{abstract}
We propose a multiscale model predictive control (MPC) framework for stationary battery systems that exploits high-fidelity models to trade-off short-term economic incentives provided by energy and frequency regulation (FR) markets and long-term degradation effects. We find that the MPC framework can drastically reduce long-term degradation while properly responding to FR and energy market signals (compared to MPC formulations that use low-fidelity models). Our results also provide evidence that sophisticated battery models can be embedded within closed-loop MPC simulations by using modern nonlinear programming solvers (we provide an efficient and easy-to-use implementation in {\tt Julia}). We use insights obtained with our simulations to design a low-complexity MPC formulation that matches the behavior obtained with high-fidelity models. This is done by designing a suitable terminal penalty term that implicitly captures long-term degradation. The results suggest that complex degradation behavior can be accounted for in low-complexity MPC formulations by properly designing the cost function. \mynote{We believe that our proof-of-concept results can be of industrial relevance, as battery vendors are seeking to participate in fast-changing electricity markets while maintaining asset integrity}.
\end{abstract}

\section{Introduction}
Batteries are flexible assets that can help modulate power grid loads at multiple timescales. A particular source of flexibility that is becoming increasingly valuable to the power grid is frequency regulation (FR) \cite{zavalastranded}. Under an FR market, the power grid remunerates a battery for providing a flexibility band that is used to modulate loads at time resolutions of seconds. From the battery perspective, determining an optimal amount of FR capacity to be offered in the market is a non-trivial task. Specifically, dynamics of FR signals can be rather aggressive and significantly deteriorate the battery life (capacity fade). Moreover, the battery needs to determine how to best use stored energy and when to buy power to replenish the battery. This involves a complex multiscale decision-making problem in which the battery must balance short-term revenue with long-term asset degradation.  In lithium-ion battery systems, one of the main reaction mechanisms of capacity fade is caused by irreversible side reactions that occur at the boundary of the electrode and electrolyte. \mynote{These side reactions create a layer known as the solid electrolyte interphase (SEI) layer \cite{peled1979electrochemical}. The growth of the SEI layer causes growth of resistance and loss of active lithium material \cite{balbuena2004lithium,verma2010review}.} Capacity fade is thus a critical consideration in battery market participation and provision of flexibility.  

Optimal participation strategies for battery systems in FR markets \cite{He2016,MohsenianRad2016,xu2018optimal} and in demand charge mitigation \cite{Lucas2016,Sebastian2016,kumar2018stochastic} have been explored in the literature. \mynote{A common limitation of these studies is that they use empirical/low-fidelity (e.g., equivalent-circuit \cite{schmalstieg2014holistic}) models. Empirical models are not able to accurately capture dynamic behavior and are formulated based on a limited number of experimental conditions (they have limited generalizability). For example, the equivalent-circuit reported in \cite{schmalstieg2014holistic} fits the data reported in the same paper fairly well (root-mean-square error of 3.14\%) but does not fit the data reported in a different paper (root-mean-square error of 12.23\%) \cite{reniers2018improving}.} More importantly, the lifetime prediction of empirical models often rely on low-fidelity representations of battery degradation (e.g., cycle counting) or impose conservative constraints that try to indirectly prevent degradation. As a result, such formulations cannot accurately capture safety constraints (e.g., maximum voltage and lifetime) and can make suboptimal market participation decisions. 

Despite the disadvantages of low-fidelity battery models, the use of such models has been motivated by the computational complexity of high-fidelity (physics-based) models, which comprise sets of highly nonlinear differential equations. A simple approach to capture capacity fade in control and optimization formulations is to consider this a function of cumulative energy \cite{walawalkar2007economics, moghaddam2018optimal}. More sophisticated approaches take into account factors such as depth of discharge (DOD) \cite{dogger2011characterization} and state of charge (SOC) \cite{moura2011optimal}. The two-step approach proposed in \cite{perez2016effect} uses an empirical model that ignores capacity fade for decision-making and employs a detailed physics-based model to determine the capacity fade incurred under such a decision (a posteriori).  \mynote{Aging-aware empirical models have also been used to address multiple case studies in power grids \cite{mosca2018battery, maheshwari2020optimizing}}.

\mynote{A few researchers have used high-fidelity models to capture capacity fade. The work in \cite{reniers2018improving} compared three battery models to perform price arbitrage, and noted that the high-fidelity model, such as the single particle models used in this paper, improved the revenue and capacity prediction error substantially compared with two aging-aware empirical models. The authors in \cite{perez2016optimal} used a high-fidelity model to determine the tradeoffs between capacity fade and charge time. The authors in \cite{pathak2017analyzing} used a high-fidelity model coupled with an optimal model-based controller for health-aware battery charging. Other researchers have used highly nonlinear physics-based battery models in renewable grid systems \cite{weisshar2015model, lee2017direct, bonkile2019power}. To the best of our knowledge, however, approaches that directly embed high-fidelity models in {\em frequency regulation market participation strategies} have not been reported in the literature.} We attribute this not only to the computational complexity of physics-based models but also to the inherent multiscale nature of the battery management problem. Specifically, battery management systems must capture long-term capacity fade effects and short-term fluctuations of electricity prices and FR signals. 

This paper presents a model predictive control (MPC) framework to simultaneously optimize FR market participation while mitigating capacity fade using high-fidelity battery models. Physics-based models can provide accurate predictions of internal states of the battery system and such states can be linked to degradation/capacity fade, thereby predicting the lifetime of the battery under dynamic operating conditions in FR markets \cite{lawder2014battery,lee2017direct}. Our framework solves a short-term (1 hour) optimization problem at high time resolution (2 seconds) and uses a terminal cost penalty on capacity fade to capture long-term effects. We conduct extensive closed-loop simulations and find that the MPC formulation provides substantial improvements in economic potential and capacity fade over formulations that use low-fidelity models. This is the result of having direct control over internal battery states. Our simulations are enabled by the use of computationally efficient discretization schemes and sparse nonlinear programming solvers. We use the knowledge gained with high-fidelity MPC simulations to design a computationally more tractable reformulation of the MPC problem that does not require a high-fidelity model and we show that this formulation provides satisfactory economic and degradation performance.  We also provide an efficient and easy-to-use {\tt Julia} implementation of the MPC framework. 

\section{Single Particle Battery Model}\label{sec:battery_model}
This section describes a model for a lithium-ion battery that will be used as the core of the proposed MPC framework. The single particle (SP) model \cite{haran1998determination, santhanagopalan2006review} used in this paper describes the micro-scale of the battery system. The SP model is obtained by neglecting the solution phase of the cell and assuming constant and uniform concentration of the electrolyte and potential in the solution phase. Although the SP model is less accurate than other full order electrochemical models (e.g., Doyle-Fuller-Newman model) when the C-rate is high, but it is also computationally more tractable. \mynote{Here, C-rate is a measure of charge and discharge rate; for example, a C-rate of 1C means that the discharge current will discharge the entire battery in an hour, while a C-rate of 2C means that the discharge current will discharge the entire battery in half an hour}. This paper uses the parameters reported in \cite{forman2012genetic}, \mynote{which have been estimated using data collected from extensive voltage/current cycling tests of $LiFePO_4$ cells. The model used assumes that the temperature of the battery is kept constant; this is reasonable to assume in stationary battery systems (because such systems have efficient cooling capabilities).}. 

In the SP model, both electrodes are assumed to be made of uniform spherical particles with a radius $R_j$,  where the subscript $j \in \{n, p\}$ represents the negative and positive electrodes, respectively. The diffusion of the lithium ions within the particles is described by Fick's second law in spherical coordinates:
\begin{align}\label{Eq:diffusion}
\frac{\partial c_j}{\partial t}  = \frac{1}{r^2} \frac{\partial}{\partial r} \left (r^2 D_j \frac{\partial c_j}{\partial r} \right ),   \;\;\;\; j \in \{n, p\}
\end{align}
with the boundary conditions:
\begin{align}\label{Eq:diffusion_bound}
 \left ( \frac{\partial c_j}{\partial r}  \right )_{r=0} = 0, \quad \left( \frac{\partial c_j}{\partial r} \right )_{r=R_j} = - \frac{J_j}{D_j F}. 
\end{align}
where $r$ represents the radial direction coordinate, $t$ is the time dimension, $c_j$ denotes the solid phase concentration of lithium in the electrode $j$, $J_j$ denotes the average local reaction current density, $D_j$ is the diffusion coefficient of lithium in the solid phase, and $F$ is the Faraday constant. \mynote{It has been previously shown in \cite{subramanian2004boundary} that the partial differential equations \eqref{Eq:diffusion}-\eqref{Eq:diffusion_bound} can be simplified into a set of differential and algebraic equations (DAEs) by approximating the concentration profile within the sphere by a parabolic profile. The work in \cite{subramanian2001approximate} shows that the approximate model matches the exact model well.} The DAE system captures the average concentration within the particle $c^{avg}_j$ and the surface concentration $c^s_j$ as:

\begin{align}\label{Eq:diffusion_ode}
 \frac{d c^{avg}_j}{dt}  =  \frac{-3 J_j}{R_j F}\\
   c^{s}_j  = c^{avg}_j  + \frac{-J_j R_j}{5 D_j F}.
\end{align}

The local current density $J_j$ is obtained by using the Butler-Volmer (BV) kinetic expression:
\begin{align}\label{Eq:BV}
&  J_j = 2\cdot i_{0,j} \cdot \sinh \left(\frac{0.5 F }{R T} \eta_j\right) \\
& i_{0,j} = F k_j (c_{j, max}-c^s_{j})^{0.5} (c_{j}^s)^{0.5} c_e^{0.5}.
\end{align}
where $T$ denotes temperature (assumed fixed), $\eta_j$ is the local overpotential, $k_j$ is the rate constant of electrochemical reaction, $c_{j, max}$ represents the maximum concentration of lithium ions in the particles of electrode $j$, $c_e$ is the concentration of electrolyte in solution phase.  The local overpotentials driving the electrochemical reaction are given by $\eta_p = \phi_p - U_p (\theta_p)$ and $\eta_n = \phi_n - U_n (\theta_n) + R_f {I_{app}}/{S_n}$, where $\phi_j$ is the solid-phase potential of electrode $j$, $R_f$ is the resistance of the SEI film, $U_j$ is the open-circuit potential, \mynote{$S_j$ is the total electroactive surface area of electrode j}, and $I_{app}$ is the applied current passing through the cell. $I_{app}$ is defined as positive for charging process and negative for discharge process. Moreover, we have that $\theta_j = c_{j}^s/c_{j,max}, j \in \{n, p\}$. Following \cite{ramadass2004development}, we assume that capacity fade is caused by an irreversible solvent reduction reaction, which causes the formation of a resistive SEI film in the negative electrode. This mechanism results in the loss of active material and the increase of internal impedance. The authors in \cite{ramadass2004development} assume that the side reaction only occurs during charging. Following observations made in \cite{bashash2011plug, moura2013battery, moura2011optimal}, however, we assume that the side reaction occurs under both charging and discharging. We argue that this assumption is not only closer to reality but, surprisingly, also makes the model computationally more tractable (it avoids discontinuous logic that turns on/off capacity fade behavior). 

The current density for the side reaction $J_{sd}$ is governed by the Butler-Volmer kinetics, which can be simplified by assuming that the reaction is irreversible and that the change of solvent concentration is small, to obtain $J_{sd} = - i_{o,sd} \exp\left( {-F \eta_{sd}}/{R T} \right)$, where $i_{o,sd}$ denotes the exchange current density for the side reaction. Symbol $\eta_{sd}$ denotes the side reaction overpotential, which is in turn given by $\eta_{sd} = \phi_n - U_{ref} + R_f {I_{app}}/{S_n}$. Here, $U_{ref}$ represents the constant equilibrium potential of the side reaction. Symbol $R_f$ denotes the total resistance of the SEI film and is given as $R_f = R_{SEI} + {\delta_f}/{\kappa_{sd}}$, where $R_{SEI}$ denotes the initial film resistance, $\kappa_{sd}$ denotes the conductivity of the film and $\delta_f$ denotes the film thickness. The film growth is governed by the differential equation:
\begin{align}\label{Eq:delta_film}
\frac{ d \delta_f} {d t} =  \frac{- J_{sd} M_{sd}}{\rho_{sd} F}
\end{align}
where $M_{sd}$ denotes the molecular weight of the side product and $\rho_{sd}$ represents the density of the side product.  The {\em rate of capacity fade} $C_r$ is a function of the maximum capacity of the cell $Q_{max}$ and is given by $C_r = \frac{J_{sd} S_n} {Q_{max}}$. The cumulative capacity fade is the integral of the fade rate and given by $C_f = \int C_r dt$.  The voltage $V$, current $I_{app}$, power $P$, and energy $E$ are computed from:
\begin{align}\label{Eq:voltage}
V&= \phi_p - \phi_n\\
J_p &= \frac{I_{app}}{S_p}\\
J_n+J_{sd} &= \frac{-I_{app}}{S_n} \\
P &=  I_{app} V / 10^6 \\
 E &= \frac{c^{avg}_n}{c_{n,{max}}}E_{max},
\end{align}
where $E_{max}$ is the battery capacity.  As can be seen, the battery model comprises a complex set of highly nonlinear differential and algebraic equations.  

\section{Multiscale Market Participation Problem}
We begin by describing the decision-making setting under which the battery is operated and we then describe the MPC formulation to automate market participation decisions. 

\subsection{Decision-Making Setting}

\begin{figure}[!htb]
\centering
\includegraphics[width=0.7\textwidth]{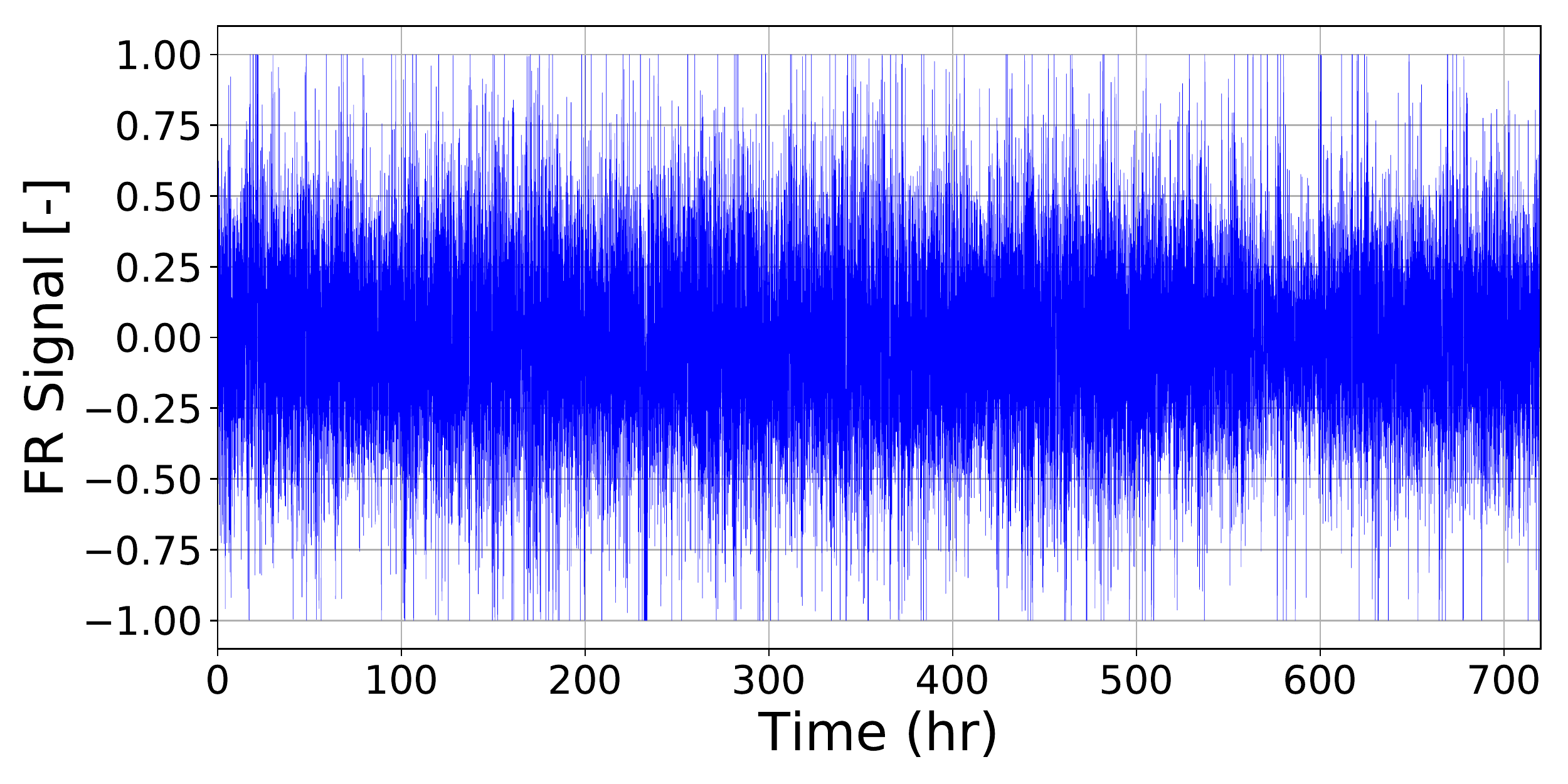} 
\includegraphics[width=0.7\textwidth]{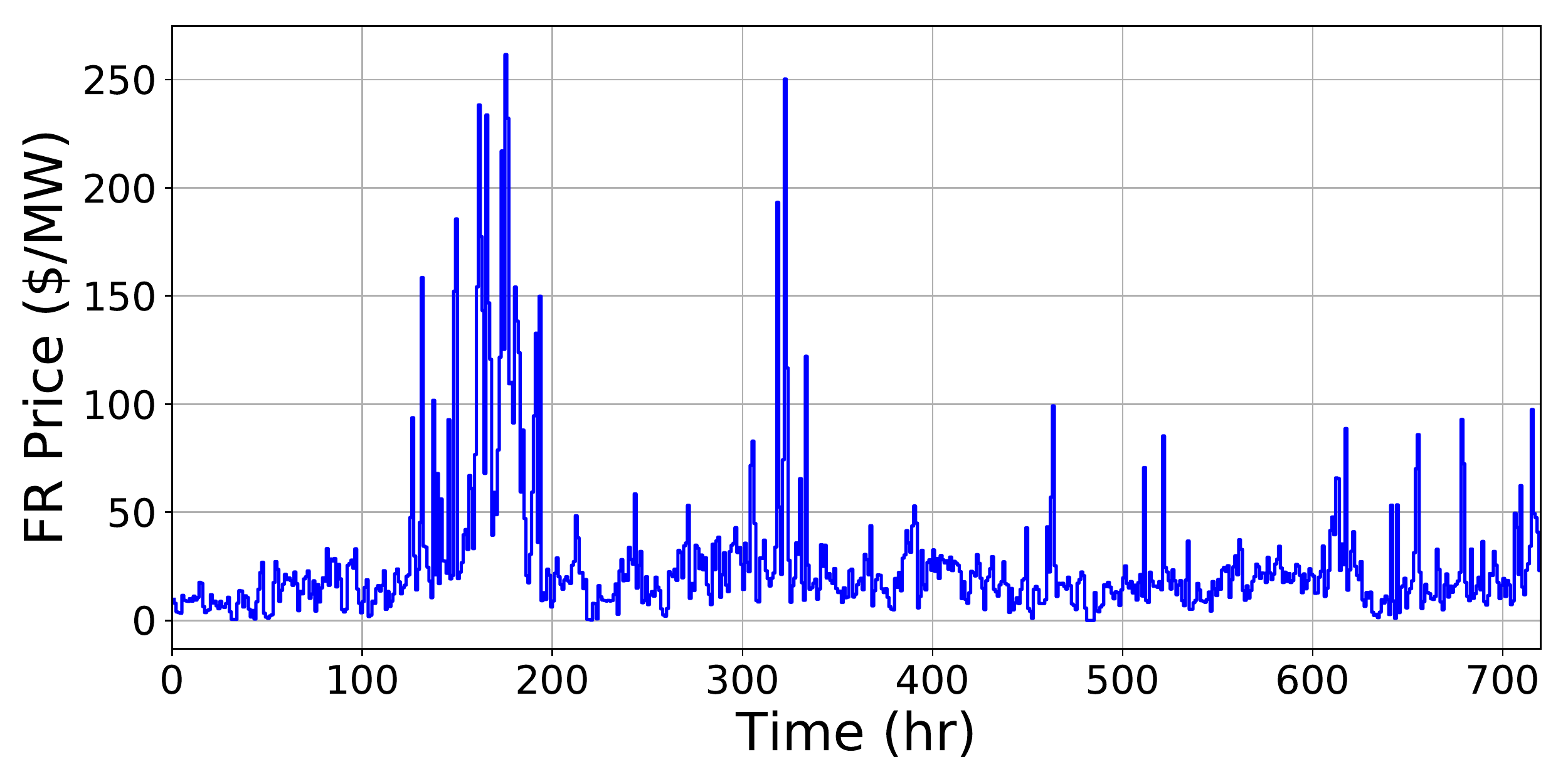}
\includegraphics[width=0.7\textwidth]{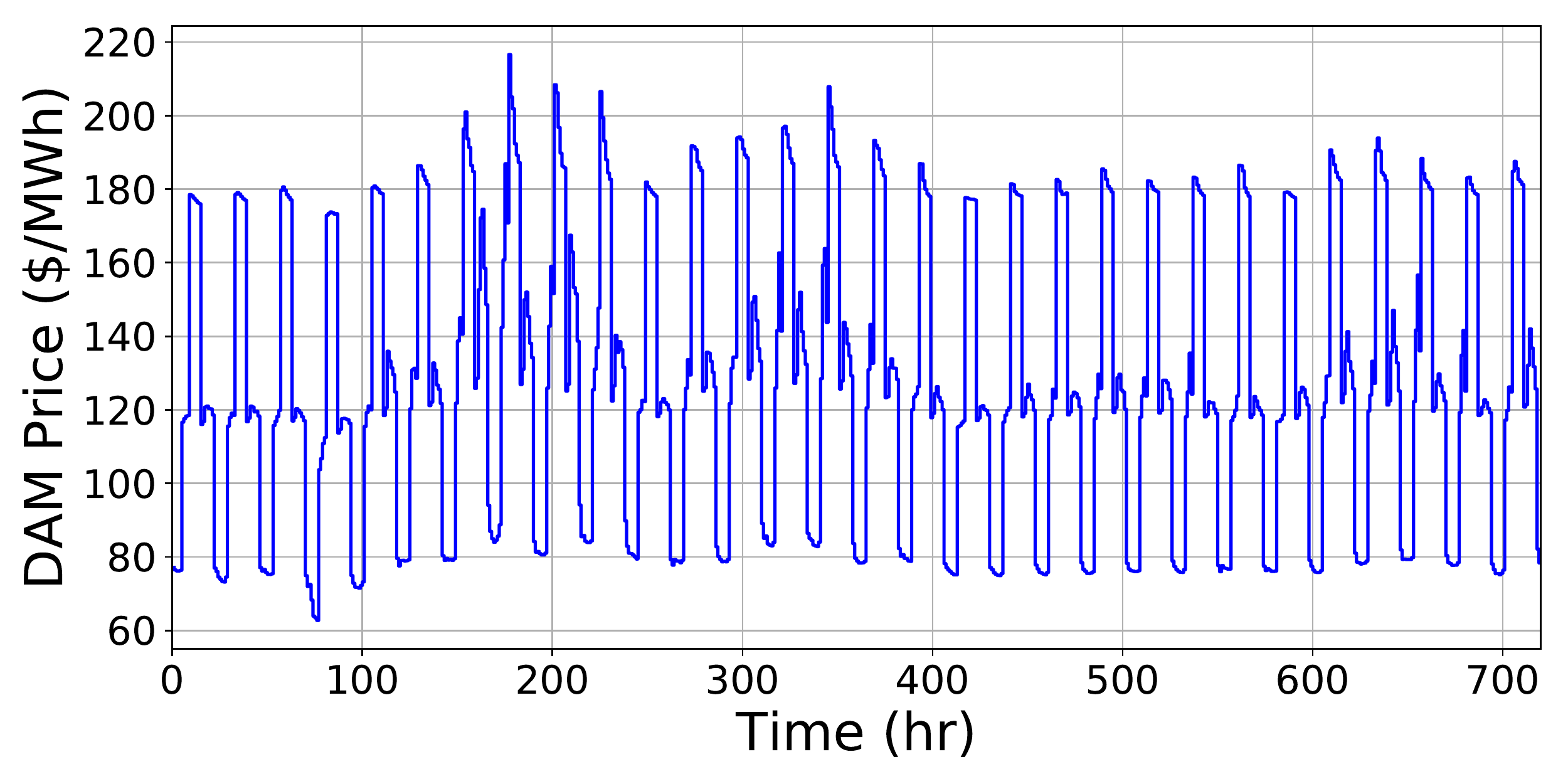}
\vspace{0.0in}\caption{PJM data for a month for FR signal (top), FR capacity price (middle) and day-ahead price (bottom).}
\label{fig:data}
\vspace{-0.05in}
\end{figure}

The battery seeks to determine optimal market participation strategies in energy and FR markets that are operated by an independent system operator (ISO), while simultaneously mitigating battery degradation. This work focuses on the setting provided by PJM Interconnection. We use real price and FR signal data from PJM to conduct our study. The various cost and revenue components that are considered are:
\begin{itemize}
\item {\em Frequency Regulation Capacity (hourly)}: 
The battery needs to decide the committed FR capacity band for the next immediate hour. The ISO can request the battery to dispatch a fraction of the committed capacity based on the grid requirements in real-time (every {\em two seconds} in PJM). The real-time FR signal from the ISO has a bounded range of [-1,+1] (see Figure \ref{fig:data}). The ISO compensates the battery for providing an operational band based on time-varying market FR capacity prices (updated every hour). In the studied setting, we ignore performance-based compensation from FR markets \cite{Chen2015,Sadeghi2017}. This assumption is motivated by recent work, which has found that FR capacity payments are significantly more lucrative than FR mileage payments \cite{dowling2017economic}. The FR capacity band provided is updated every hour and remains constant over the entire hour. 

\item {\em Power Purchase (hourly)}: Power can be purchased from the day-ahead energy market (DAM) to recharge the battery. The optimal purchase timing is driven by a time-varying market price (updated every hour). The amount of power purchased can be updated every hour and remains constant over the hour. 

\item {\em Load (hourly)}: Power can be withdrawn by an adjustable load to help maintain the amount of energy remaining in the battery. We assume that the load can be updated every hour and remains constant over the entire hour. The revenue associated with this energy load is zero. 

\end{itemize}

\subsection{High-Fidelity MPC Formulation}

The battery management problem that is tackled in this work is {\em multiscale} in nature because it must capture {\em hourly} variations in energy and FR price signals, {\em second-by-second} variations of FR signals, and {\em long-term} battery degradation (spanning days to years). Under the proposed MPC framework, that we call high-fidelity MPC (HF-MPC), an optimization problem is solved at every hour $t$ over the prediction horizon $\mathcal{N}_t:=\{t+1,t+2,...,t+N\}$. Here, $N$ is the length of the prediction horizon. Since the FR signal is updated every two seconds, each hour is discretized using $S$=1,800 time steps and we define the time interval set $\mathcal{S}:=\{1,...,S\}$.  

The parameters of the MPC formulation are: $\pi^e_k\in \mathbb{R}$ denotes the electricity price at $k$th hour [\$/MWh], $\pi^f_k\in \mathbb{R}_+$ is the FR capacity price at $k$th hour [\$/MW],  $\alpha_{k,s} \in [-1,1]$ is the fraction of FR capacity requested by the ISO at $k$-th hour and $s$-th step [-] (if $\alpha_{k,s}>0$, the ISO sends power while if $\alpha_{k,s}<0$ the ISO withdraws power),  $\overline{E}\in \mathbb{R}_+$ is the battery capacity [MWh], $\overline{P}\in \mathbb{R}_+$ is the maximum charging rate [MW], and $\underline{P}\in \mathbb{R}_+$ is the maximum discharging rate  [MW]. \mynote{In this paper, we assume that the FR signal and price signals are known in advance (we consider a deterministic process). }

The variables of the MPC formulation are:  $F_k\in \mathbb{R}_+$ is the FR capacity provided at $k$-th hour [MW], $O_k\in \mathbb{R}_+$ is the power purchased from the day-ahead-market at $k$-th hour [MW], $L_k\in \mathbb{R}_+$ is the committed power to load at $k$-th hour  [MW], $P_{k,s}\in \mathbb{R}$ is the net battery charge/discharge rate at $k$-th hour and $s$-th step [MW] ($P_{k,s}>0$ the battery is being charged and if $P_{k,s}<0$ the battery is being discharged),  $x_{k,s}\in \mathbb{R}_+$ are the state variables of battery at $k$-th hour and $s$-th step ($c^{avg}_j$, and $\delta_f$), $E_{k,s}\in \mathbb{R}_+$ is the remaining energy in the battery [MWh], $C^f_{k,s}\in \mathbb{R}_+$ is the capacity fade [-],  and $V_{k,s}\in \mathbb{R}_+$ [V] is the voltage. \mynote{This paper assumes that the states $x_{k,s}$ can be estimated accurately from voltage and current data. The estimation of states is a separate (and challenging) research topic and is not explored in this study}.

All quantities with a single subindex $k$ are held constant over the time interval $[(k-1),k]$ and all quantities with subindices $k,s$ are held constant over the interval $[k-1+(s-1)/S, k-1+s/S]$. The FR capacity $F_k$ represents a {\em symmetric band} and the actual FR power requested by the ISO is $-\alpha_{k,s}F_k$. 

\subsubsection{Objective Function}\label{subsubsec:obj_uncert}
The objective of the MPC problem is to maximize profit (considering the revenue from FR participation and the energy cost) while penalizing capacity fade over the horizon $\mathcal{N}_t$:
\begin{align}\label{Eq:objective}
 &\sum\limits_{k\in \mathcal{N}_t} {\pi}^f_{k} F_{k}  - \sum\limits_{k\in \mathcal{N}_t} {\pi}^e_{k} O_{k}  - {\pi}^{C^f} (C^f_{t+N,S} - C^f_{t+1,1})
\end{align}
 The first term is the revenue obtained from the provision of FR {\em capacity}, the second term is the cost of purchasing power from the day-ahead market, and the third term is the capacity fade penalty. The parameter ${\pi}^{C^f}$ is the penalty parameter, which estimates the long-term value of capacity fade.  

\subsubsection{Constraints}\label{subsubsec:const_uncert}
The battery model introduced in Section \ref{sec:battery_model} is discretized using backward Euler scheme in time and can be expressed in the following compact form:
\begin{align}
x_{k,s+1} = \varphi_1 ( x_{k,s}, P_{k,s}), \;\; k \in \mathcal{N}_t, s \in  \mathcal{S} \\
(E_{k,s}, C^f_{k,s},V_{k,s}) = \varphi_2(x_{k,s}, P_{k,s}), \;\; k \in \mathcal{N}_t, s \in  \mathcal{S}
\end{align}
The net charged/discharged battery power equals the amount of power sent to the battery due to FR participation plus the amount of energy ordered minus the amount of load 
\begin{align}
P_{k,s} = \alpha_{k,s} F_k + O_k - L_k, \;\; k \in \mathcal{N}_t, s \in  \mathcal{S}.
\end{align}
The use of high-fidelity model enables us to directly impose safety constraints on internal states such as the voltage $\underline{V} \leq V_{k,s} \leq \overline{V}, \, k \in \mathcal{N}_t, s \in  \mathcal{S}$. 	
Due to capacity fade, the amount of remaining capacity is given by $(1-C^f_{k,s})\overline{E}$. The following constraint is used to ensure that the stored energy is within the remaining capacity: $\tau_l (1-C^f_{k,s}) \overline{E}  \leq E_{k,s} \leq  \tau_u (1-C^f_{k,s}) \overline{E}, \;\; k \in \mathcal{N}_t, s \in  \mathcal{S}$.  Parameters $\tau_l$ and $\tau_u$ impose a safety margin to prevent over-charge and over-discharge. We consider a terminal constraint on the remaining energy at the end of the prediction horizon: $\eta_l (1-C^f_{k,s}) \overline{E}  \leq E_{t+N,S} \leq  \eta_u (1-C^f_{k,s}) \overline{E}$. \mynote{The terminal constraint is enforced to prepare the battery for market participation in the next horizon.} We also impose simple logical bounds on variables: $ 
-\underline{P}  \leq P_{k,s}  \leq \overline{P}$ and $0  \leq F_k  \leq \overline{P},\, k \in \mathcal{N}_t, s \in  \mathcal{S}$.

\subsubsection{Implementation}\label{subsubsec:implement}
The problem at time $t$ uses data over the prediction horizon $\alpha_{\mathcal{N}_t}$, $\pi_{\mathcal{N}_t}^f$, and $\pi_{\mathcal{N}_t}^e$. The problem solved at time $t$ is denoted as $\mathcal{P}_t(\alpha_{\mathcal{N}_t},\pi_{\mathcal{N}_t}^f,\pi_{\mathcal{N}_t}^e,x_{t,S})$. For convenience, we simplify notation and state the problem as $\mathcal{P}_t(x_{t,S})$. The solution of this problem yields the optimal commitments $P_{\mathcal{N}_t}$, $O_{\mathcal{N}_t}$ and $L_{\mathcal{N}_t}$. Only the commitments of the next hour $P_{t+1}$, $O_{t+1}$ and $L_{t+1}$ are implemented, the horizon is shifted by one hour and the problem is solved again. The MPC scheme runs over a two-year period $\mathcal{Y} = [1,...,Y]$ (with $Y$=17,520  hours) or until the battery end of life (EOL), which is defined as the elapsed time before capacity fade reaches 20\%.   The implementation is summarized as follows: 
\begin{itemize}
\item START at $t = 0$ with $x_0$ corresponding to a new half-charged battery. 
\item SOLVE $\mathcal{P}_t(x_{t,S})$ by using the data $\alpha_{\mathcal{N}_t}$, $\pi_{\mathcal{N}_t}^f$, and $\pi_{\mathcal{N}_t}^e$ to obtain commitments $F_{t+1}$, $O_{t+1}$ and $L_{t+1}$.  \label{step2}
\item INJECT decisions over $(t,t+1)$. COMPUTE the net battery charge/discharge rate $P_{t+1,s} = \alpha_{t+1,s} F^{b}_{t+1} + O_{t+1} - L_{t+1}$. With $P_{t+1,s}$ and $x_t$, simulate the battery dynamics using a high-fidelity DAE simulator to obtain the UPDATED state $x_{t+1}$ and $C^f_{t+1,S}$. 
\item If $C^f_{t+1,S}\geq  0.2$, set EOL $= t$, BREAK. 
\item Set $t \leftarrow t+1$, RETURN to Step \ref{step2}
\end{itemize}

We compare different battery management strategies based on the total amount of profit earned before the end of life: 
\begin{align}  	
\Phi = \sum\limits_{k=1}^{EOL} {\pi}^f_{k} F^b_{k}  - {\pi}^e_{k} O_{k}.
\end{align}
Every closed-loop MPC simulation requires the solution of {\em tens of thousands of nonlinear programs} (which embed the battery model).  This task is computationally expensive and requires efficient solvers. 

\subsection{Low-Fidelity MPC Strategy}\label{subsubsec:simplified}
To establish a comparison, we also consider a simplified MPC formulation based on a low-fidelity battery model. We call this strategy low-fidelity MPC (LF-MPC). Here, the battery model is solely based on the energy balance (assuming an efficiency of 100\%). The capacity fade penalty in the objective function is removed and the battery dynamics are given by $E_{k,s+1} =E_{k,s}+ P_{k,s}$. Because of the simplicity of the model, the computational cost of this MPC strategy is small. The limitation of this approach is that it does not consider detailed states of the battery (e.g. current, voltage, capacity fade). Consequently, it cannot explicitly impose capacity fade and safety constraints. This type of model has been widely used in the literature \cite{kumar2018stochastic}.  A reason for this is that the problem is a linear programming problem that is easier to solve. 

\vspace{-0.1in}
\subsection{Heuristic Strategy} \label{simulation}

We considered a heuristic approach to guide battery market participation using simple decision-making logic. This is motivated by the observation that, if the prediction horizon is just one hour, the number of degrees of freedom in the problem is small (these correspond to the three commitment variables $F_{t+1}$, $O_{t+1}$, $L_{t+1}$). Consequently, it is possible to perform an exhaustive search of the space. This simulation-based approach provides sensitivity information on how profit and degradation change with the decision variables. Moreover, this approach does not require solving optimization problems. However, a large number of simulations will be needed to span the entire decision space and the method is not scalable. Consequently, we reduce the number of decision variables by using the following logic: we assume that the energy efficiency is 100\% and  we set $\eta_l=\eta_u$ in the terminal constraint. Under this assumption, the ideal remaining power at the end of prediction horizon is $E^\ast_{t+1,S}=\eta_l(1-C^f_{t,S}) \overline{E} $. If $F_{t+1}$ is known, the amount of energy from FR is $ \frac{1}{S}\sum_{s \in \mathcal{S}} \alpha_{t+1,s} F_{t+1}$. We further denote $\Delta E =E^\ast_{t+1,S} -  E_{t,S} - \frac{1}{S}\sum_{s \in \mathcal{S}} \alpha_{t+1,s} F_{t+1}$. Therefore, if $\Delta E >0$, we set $O_k=\Delta E$ and $L_k=0$ to satisfy the terminal constraint.  On the other hand, if $\Delta E <0$, we set $O_k=0$ and $L_k=-\Delta E$. In this way, the only degree of freedom to optimize for at time $t$ is the FR commitment $F_{t+1}$.  We reduce the search space for this variable by enforcing a {\em fixed FR band policy} (i.e., we search for a fixed FR band $F$). At each time step, we try the fixed FR band value by setting $F_{t+1}=F$.  If this FR band commitment results in an infeasible solution (over-charge or over-discharge), the value of the FR band is adjusted. After exploration of the entire feasible region, we determine the value that maximizes profit.  This simple logic gives insights into how FR market participation affects profit and battery degradation (it allows us to navigate inherent trade-offs). 

\begin{figure}[!htb]
\centering
\includegraphics[width=0.4\textwidth]{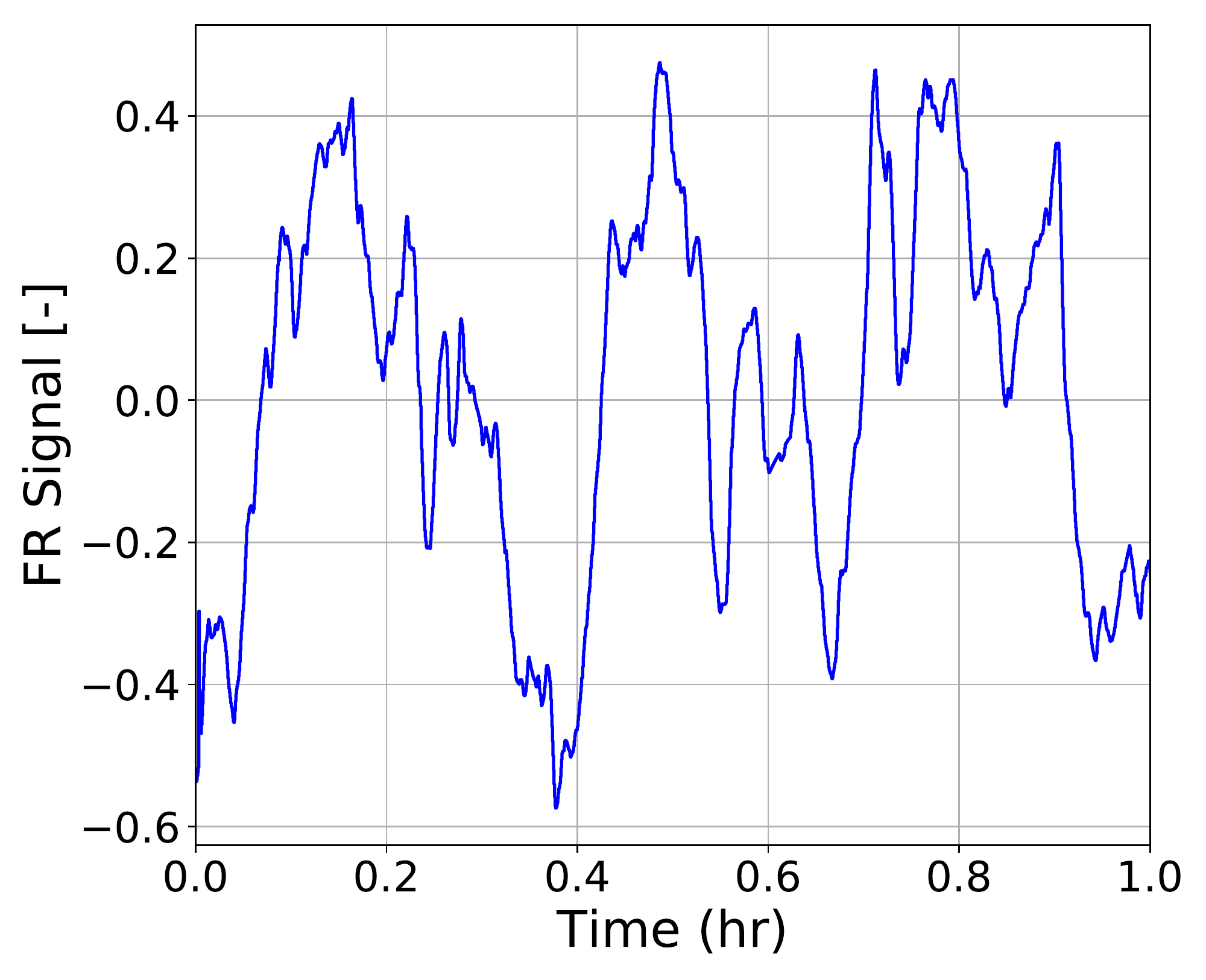} 
\label{fig:sign_hour}
\includegraphics[width=0.4\textwidth]{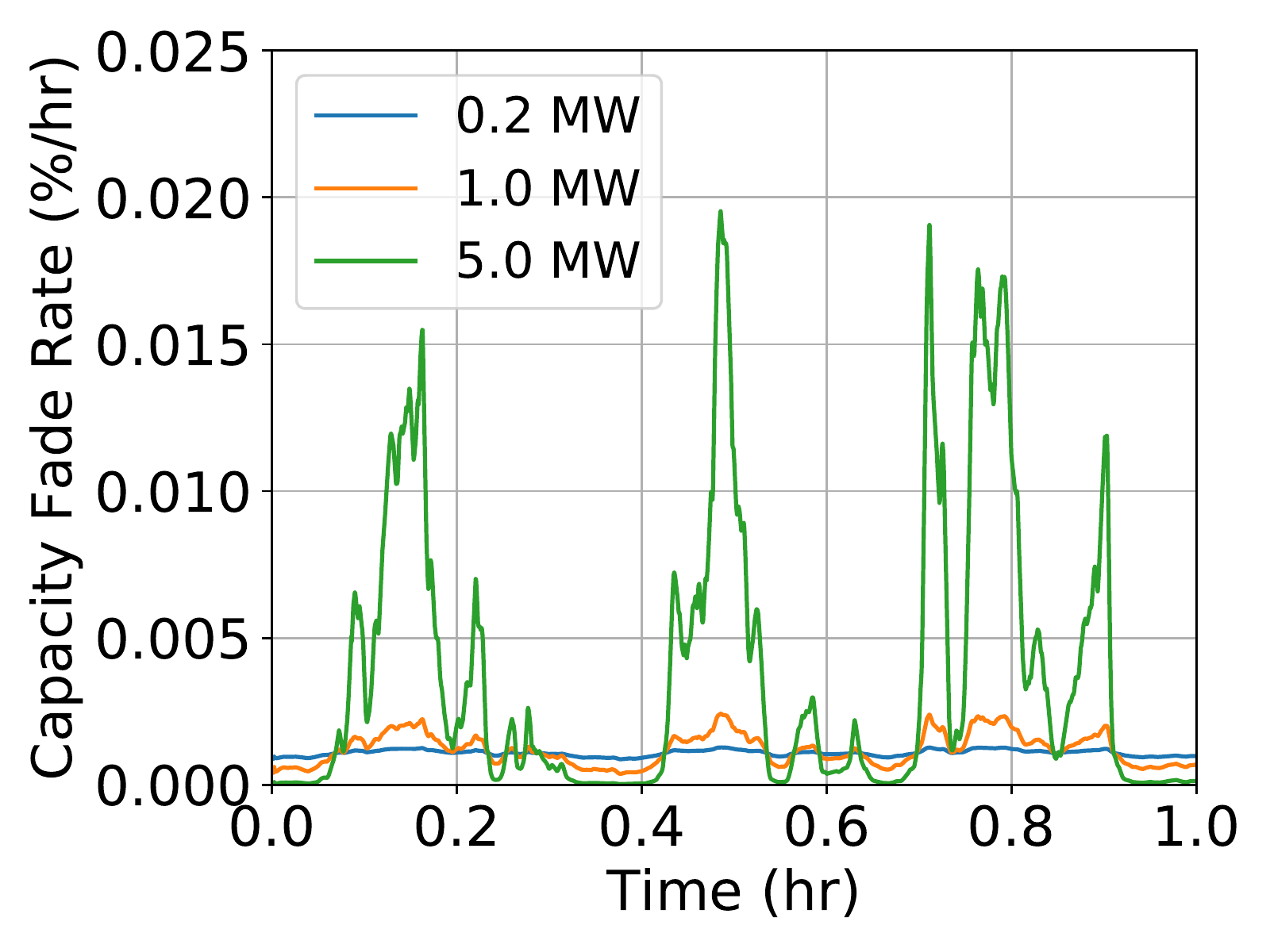}
\label{fig:fade_hour}
\includegraphics[width=0.4\textwidth]{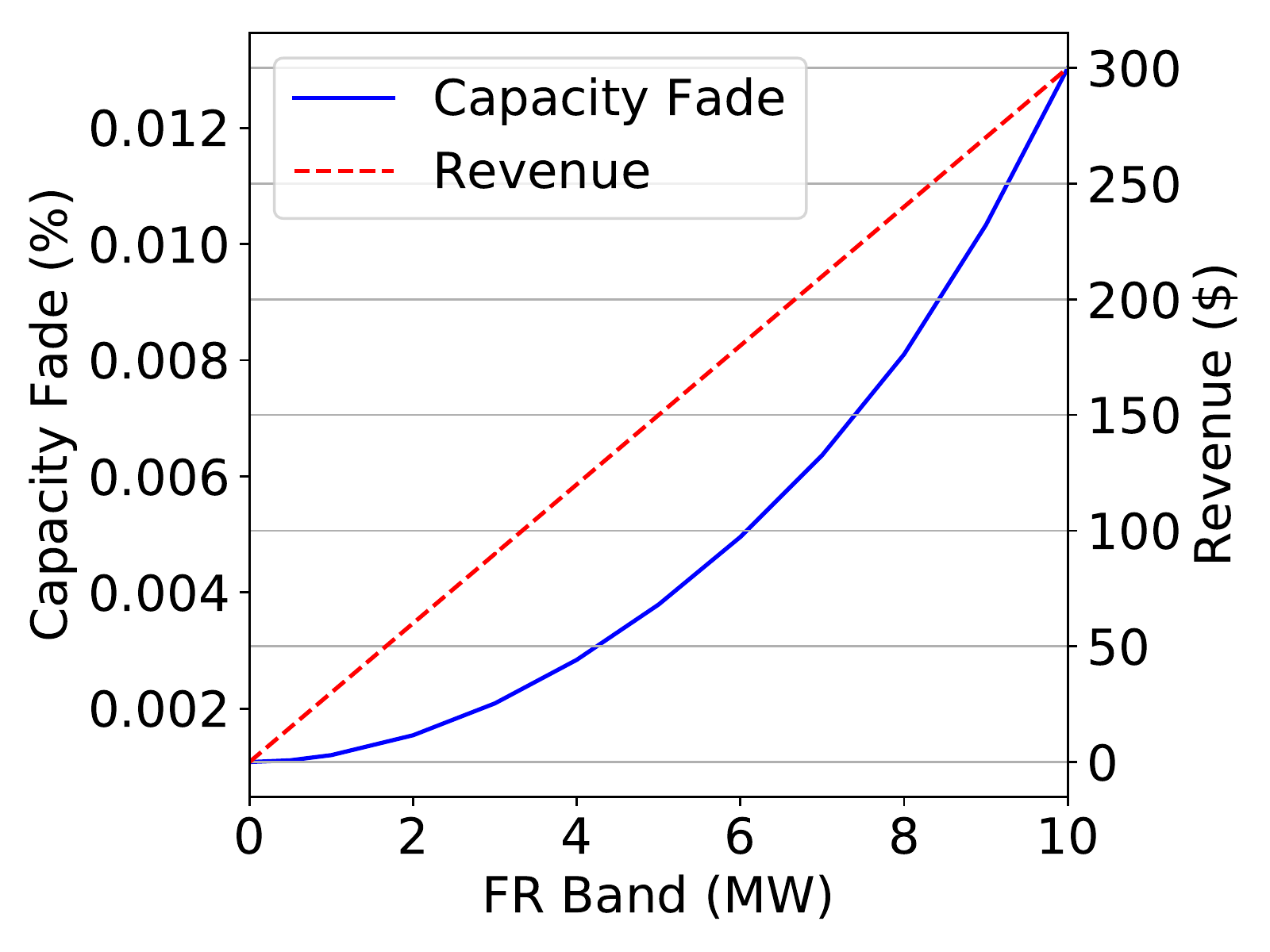}
\label{fig:fade_revenue_hour}
\includegraphics[width=0.4\textwidth]{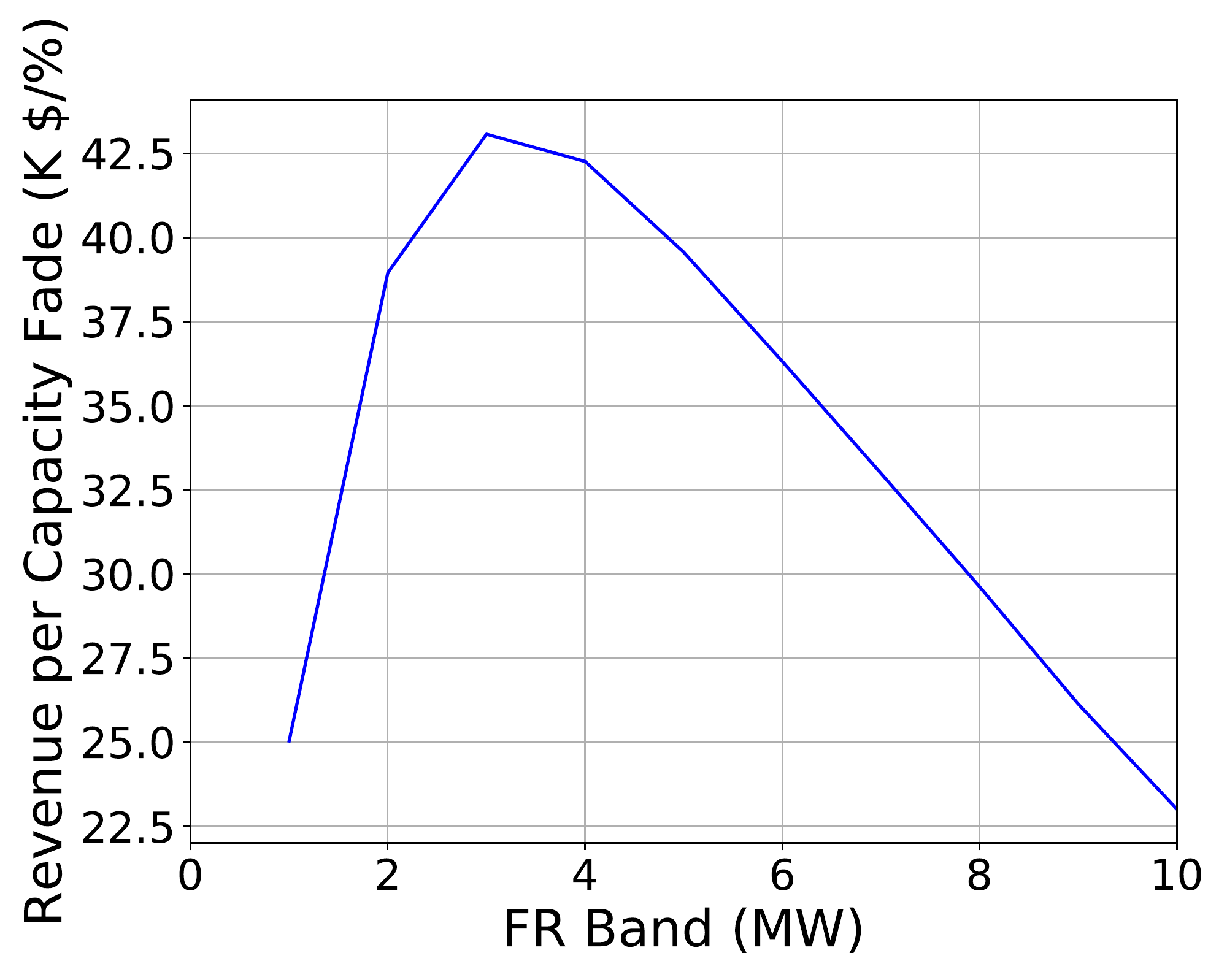}
\label{fig:fade_revenue_hour2}
\vspace{0.0in}\caption{FR signal, capacity fade rate, cumulative capacity fade and revenue, and revenue per capacity fade for one hour of simulation and different FR bands offered.}
\label{fig:hour}
\end{figure}

\section{Computational Experiments} \label{sec:experiment}

We consider a battery consisting of A123 Systems ANR26650M1 cells with lithium iron phosphate ($LiFePO_4$) cathodes. The parameters for each cell are estimated in \cite{forman2012genetic}. The number of cells is scaled so that the battery has a total capacity of 1 MWh. We assume that the maximum charging/discharging rate of the battery is 10 MW. This is a conservative assumption because the maximum continuous discharge rate of this cell is 20C, while the maximum pulse discharge rate (10 seconds) is 48C \cite{A123}. The C-rate is defined as the charge or discharge current divided by the battery capacity. We also set $\tau_l=0.1$, $\tau_u=0.9$, and $\eta_l=\eta_u=0.5$.  Because LF-MPC and the heuristic strategies cannot explicitly deal with safety constraints and a lithium iron phosphate battery has excellent safety characteristics, safety constraints are not considered in this paper. We emphasize, however, that safety constraints can be explicitly accounted for in HF-MPC. It is possible that the market commitment decisions ($F_{t+1}$, $O_{t+1}$, $P_{t+1}$) obtained with LF-MPC and the heuristic strategies are are not feasible; that is, the DAE simulator finds that the battery is over-charged or over-discharged under the computed commitments. \mynote{It is also possible that the market commitment decisions obtained with HF-MPC are  not feasible because the discretization time step used in the optimization formulation is different from that of the DAE simulator. In these cases, the value of FR band is adjusted as $F_{t+1}\leftarrow F_{t+1}-\Delta F$, where $\Delta F = 0.5$MW. Then the values of  $O_{t+1}$ and $P_{t+1}$ are re-computed.}

Historical data for one year for energy prices, FR capacity prices, and FR signals from PJM are used in  the study. Multi-year horizons are considered by replicating the historical data. All optimization problems are implemented in the modeling language {\tt JuMP}. Nonlinear programs arising in the HF-MPC strategy are solved using {\tt Ipopt}, while linear programs (LPs) arising in LF-MPC are solved using {\tt Gurobi}. All computations were performed on a multi-core computing server with Intel(R) Xeon(R) CPU E5-2698 v3  processors running at 2.30GHz.  All scripts needed to reproduce the results are available in \url{https://github.com/zavalab/JuliaBox/tree/master/Battery_FP}.  The models developed in this work have been tuned to achieve high computational efficiency (needed to deal with fast FR signal dynamics and nonlinearity). 

\subsection{Revenue vs. Capacity Fade Trade-offs} 

We first used the heuristic strategy to assess inherent trade-offs between short-term revenue and long-term battery degradation. We consider a single hour with the FR signal shown in Figure \ref{fig:hour}. Here, we can see that fast and abrupt fluctuations exist. We simulate the battery model with different committed FR band capacities while $O_k$ and $L_k$ are both set to zero. Figure \ref{fig:hour} also illustrates how capacity fade grows over time as the committed FR band increases. When the battery is charging (FR signal $\alpha>0$), increasing the committed FR band significantly increases the capacity fade rate. When the battery is discharging, capacity fade rate is relatively low. This clearly illustrates how high FR revenue can lead to faster degradation (due to the  strong fluctuations of the FR signal). We can also observe how the cumulative capacity fade and revenue change over one hour  as the committed FR band increases. In particular, we can observe that the cumulative capacity fade is positive when the FR band is zero, and it increases in a {\em nonlinear} manner. As a result, the ratio of revenue and capacity fade reaches its peak when the FR band is 3 MW. This result is important because it indicates that an optimum FR capacity indeed exists. Moreover, as we will see, the heuristic and LF-MPC strategies {\em tend to maximize FR band capacity to maximize revenue} (well above  the optimum trade-off that balances revenue and degradation). In other words, those strategies lead to aggressive market participation strategies that lead to fast degradation of the battery. We will also see that HF-MPC can correctly identify the optimal trade-off point between revenue and degradation.

\subsection{Low-Fidelity vs. High-Fidelity MPC} 

The key advantage of MPC is that it can adapt committed capacity based on market conditions and the battery internal state. To accurately capture the FR signals, we have found that it is necessary to discretize the model using timesteps of two seconds, giving rise to 1,800 time steps per hour. We consider a time horizon of one hour and 24 hours for the LF-MPC strategy and a time horizon of one hour for the HF-MPC strategy. 

The optimization problem solved at each step for the LF-MPC strategy is an LP, which contains 86,000 variables with a time horizon of 24 hours. Gurobi can solve each optimization problem in about 15 seconds. A closed-loop simulation for two years of operation requires around 14 hours of wall-clock computing time.  

The optimization problem solved in the HF-MPC strategy is a highly nonlinear NLP with 34,000 variables. On average, Ipopt requires 70 seconds to solve each optimization problem. Despite the fast solution (relative to the commitment time of one hour), a closed-loop simulation for two years requires {\em six days of wall-clock time}. Extending the horizon of HF-MPC  to 24 hours would result in an NLP with 816,000 variables. A single instance of this problem can be solved, but the closed-loop simulation would require weeks of computing time. Addressing the tractability of such formulation is an important topic of future work. Despite these limitations, we now proceed to show that dramatic improvements in economic performance and capacity fade can be achieved with HF-MPC (even with a short prediction horizon of one hour). 

\begin{figure}[!h]
\begin{center}
\includegraphics[width=0.6\textwidth]{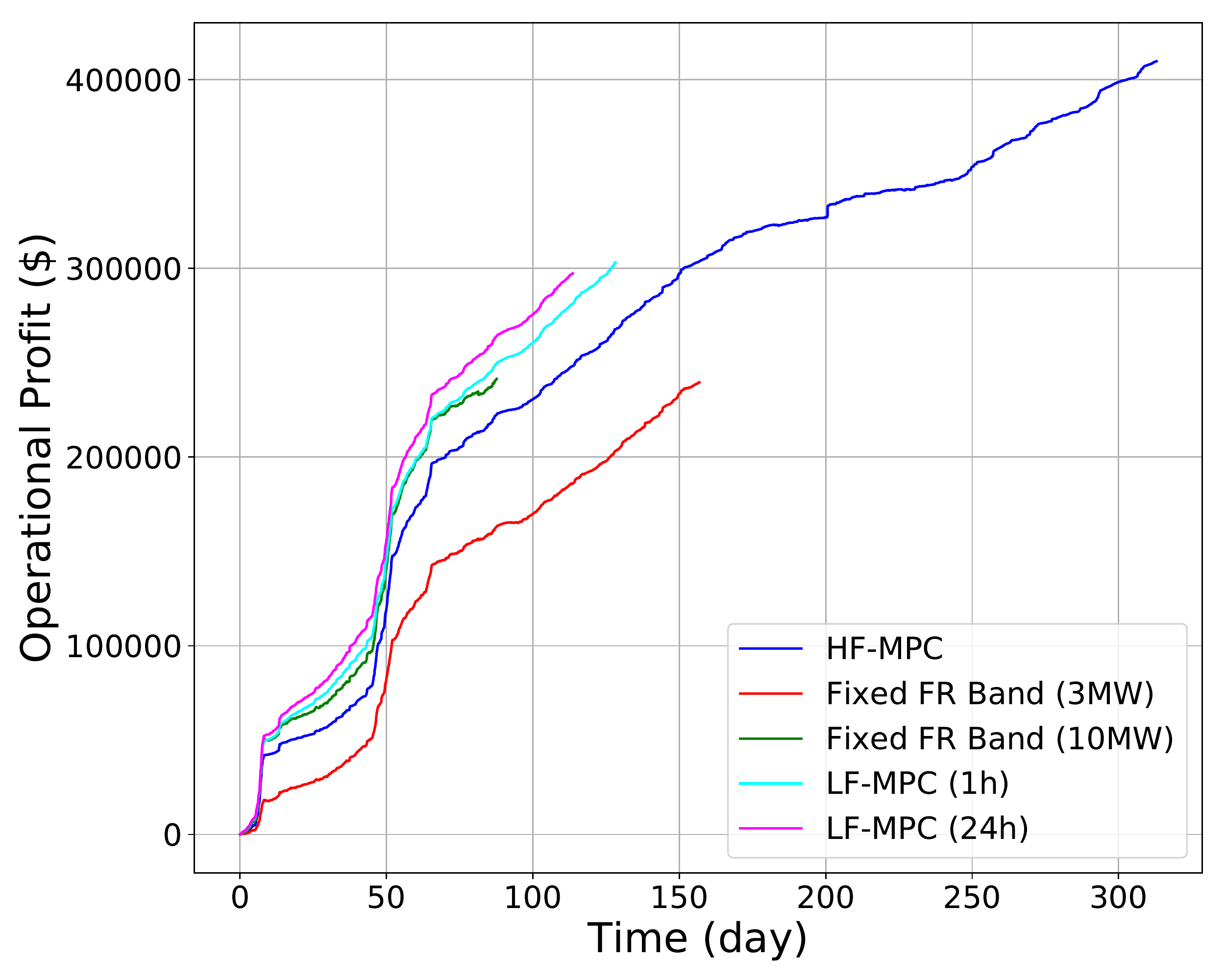}
\label{fig:operational_profit}
\includegraphics[width=0.6\textwidth]{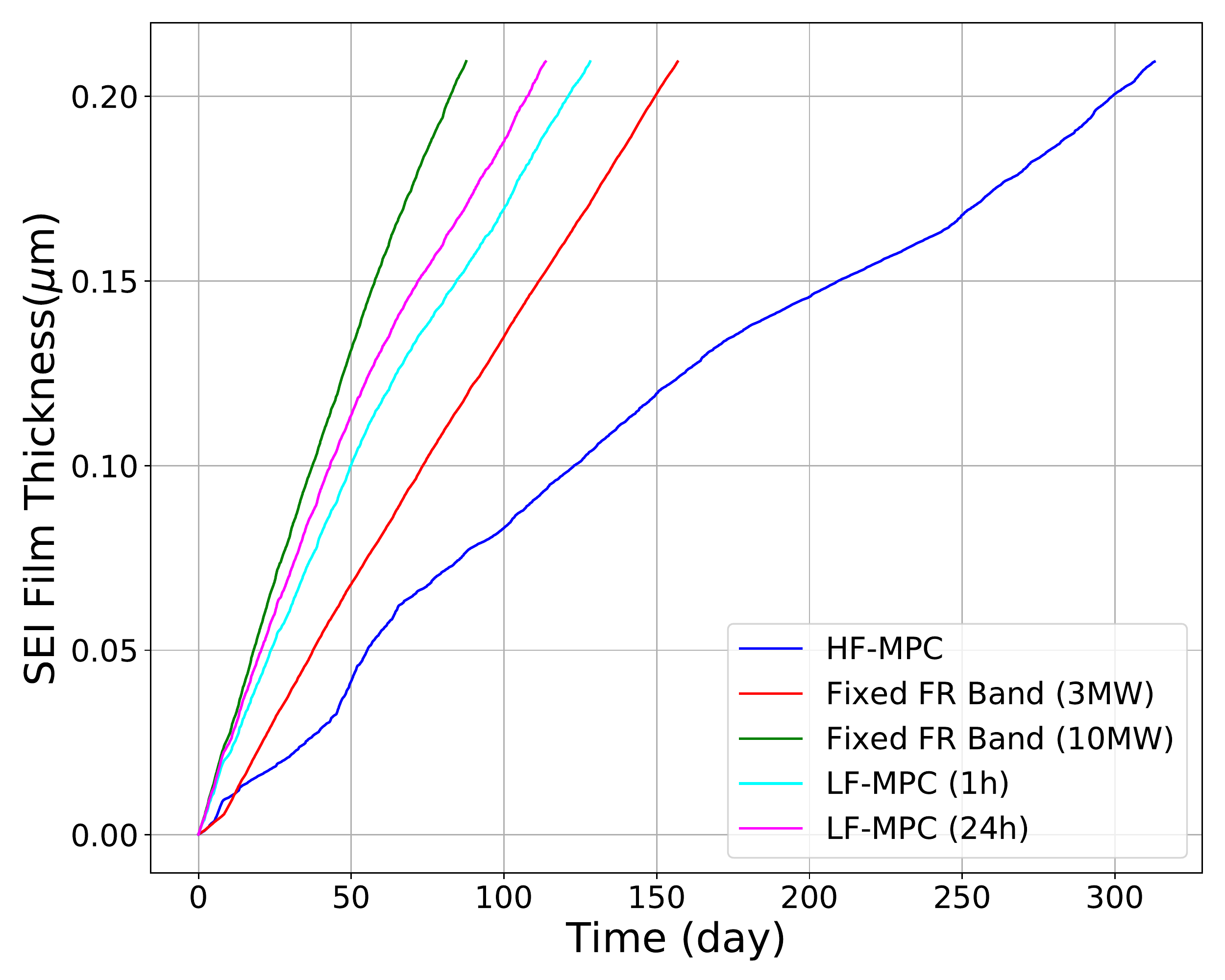}
\label{fig:film}
\caption{Profit and SEI film thickness for different strategies.}
\label{fig:comparison_year}
\end{center}
\end{figure}

Figure \ref{fig:comparison_year} shows how profit and SEI film thickness evolve over time under the different control strategies. Here, we compare the fixed FR capacity policy, the LF-MPC policy, and the HF-MPC policy. The SEI film thickness for LF-MPC grows faster than that of the fixed band policy at 10 MW (but slower than the fixed band policy at 3 MW). In contrast, HF-MPC is more strategic and offers an FR band in such a way that the SEI film thickness grows at a slower rate. As a result, the remaining capacity for this MPC policy is significantly higher. For instance, at day 100, the heuristic policy at 10 MW has reached its end of life, LF-MPC has a remaining capacity of 83\%, and HF-MPC has a remaining capacity of 92\%. The profit of HF-MPC is not significantly higher than that of other strategies for the first several months but, because the lifetime is extended significantly, the overall profit is much higher. This illustrates the ability of MPC to trade-off short-term and long-term economic performance.  

\begin{figure}[!h]
\begin{center}
\includegraphics[width=0.6\textwidth]{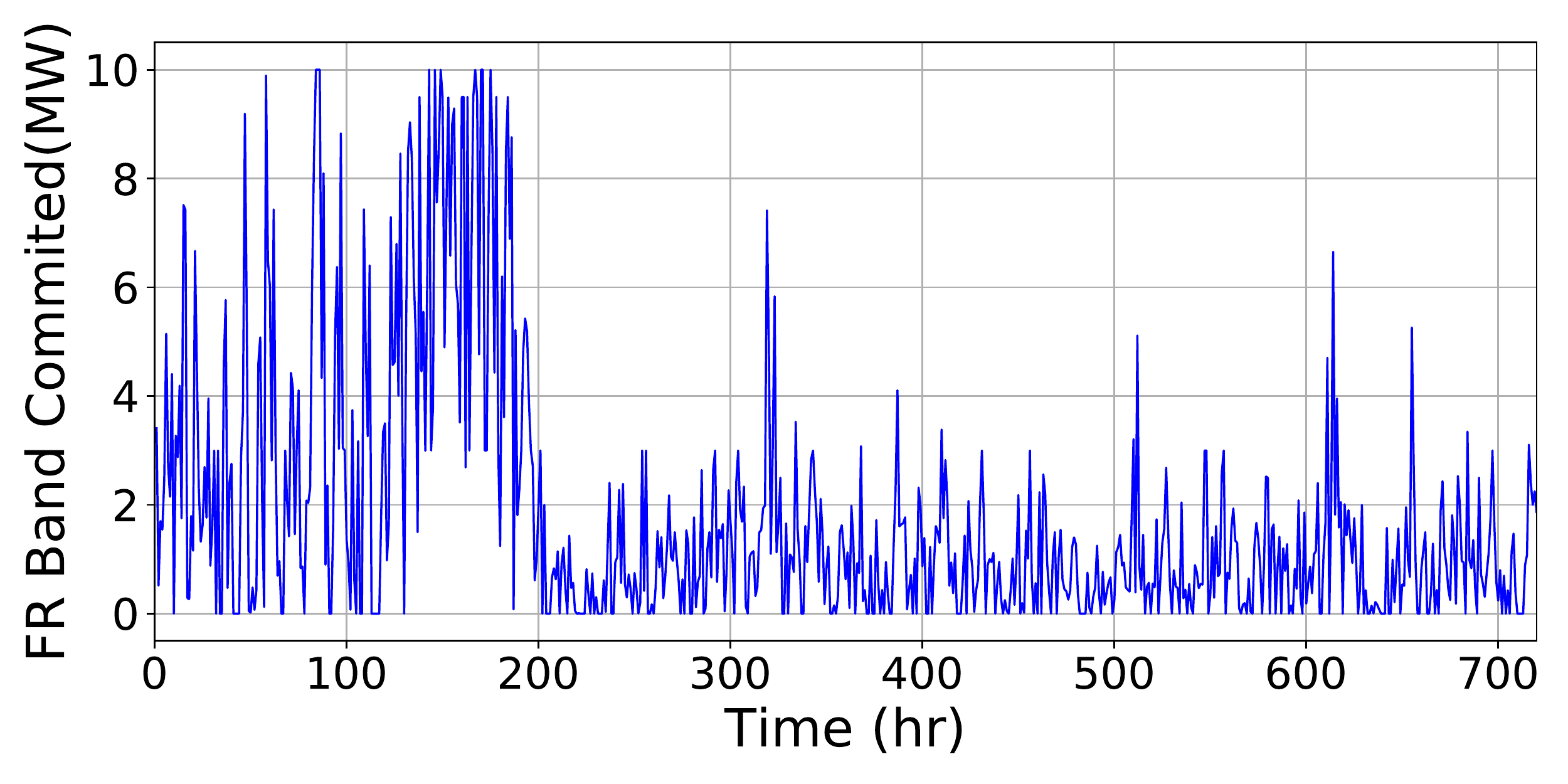}
 \\
\includegraphics[width=0.6\textwidth]{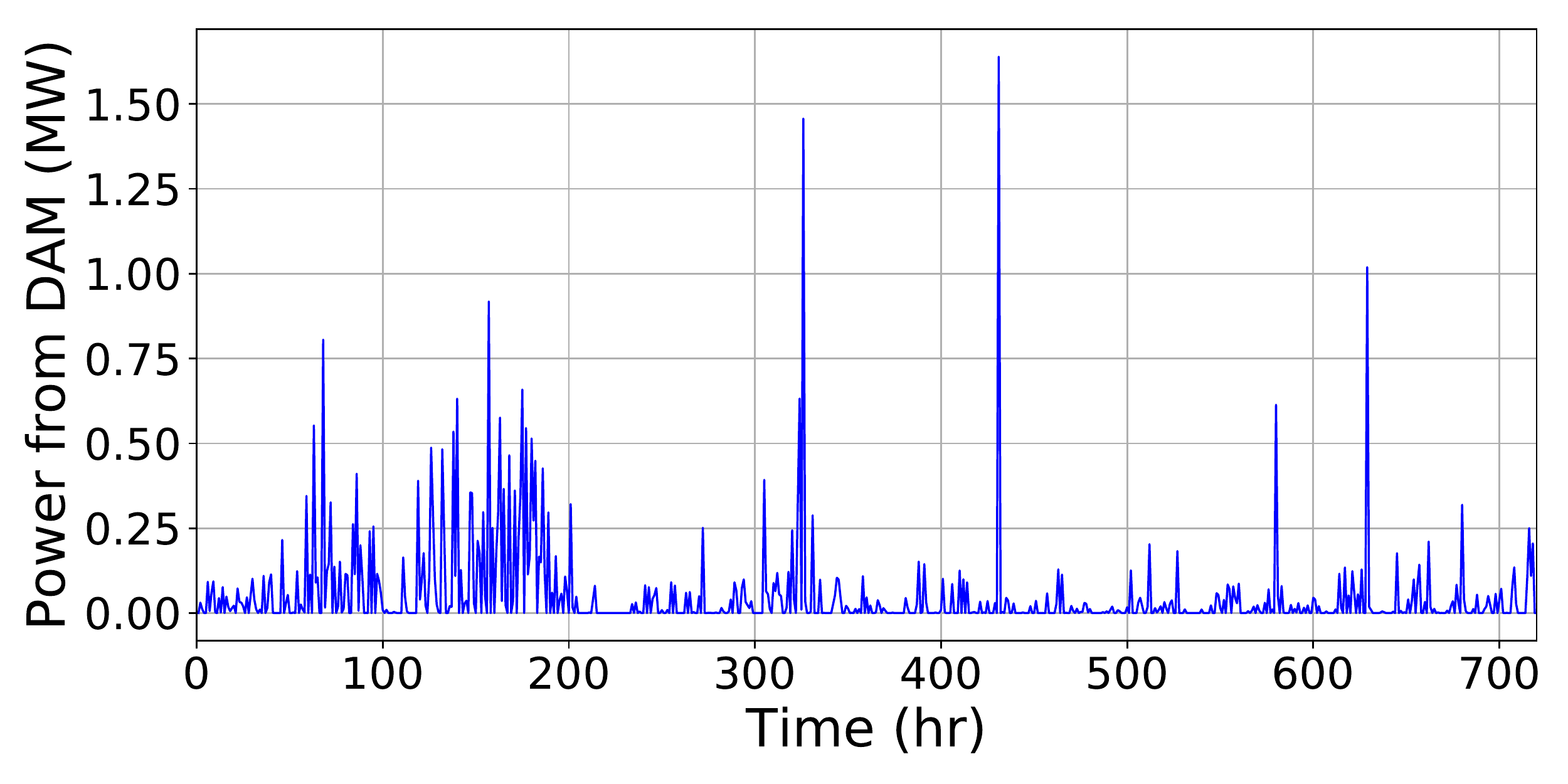}
\caption{Time profiles for commitment variables for the first month using HF-MPC.}
\label{fig:MPC_decisions_first_month}
\end{center}
\vspace{-0.2in}
\end{figure}

\begin{figure}[!h]
\begin{center}
\includegraphics[width=0.6\textwidth]{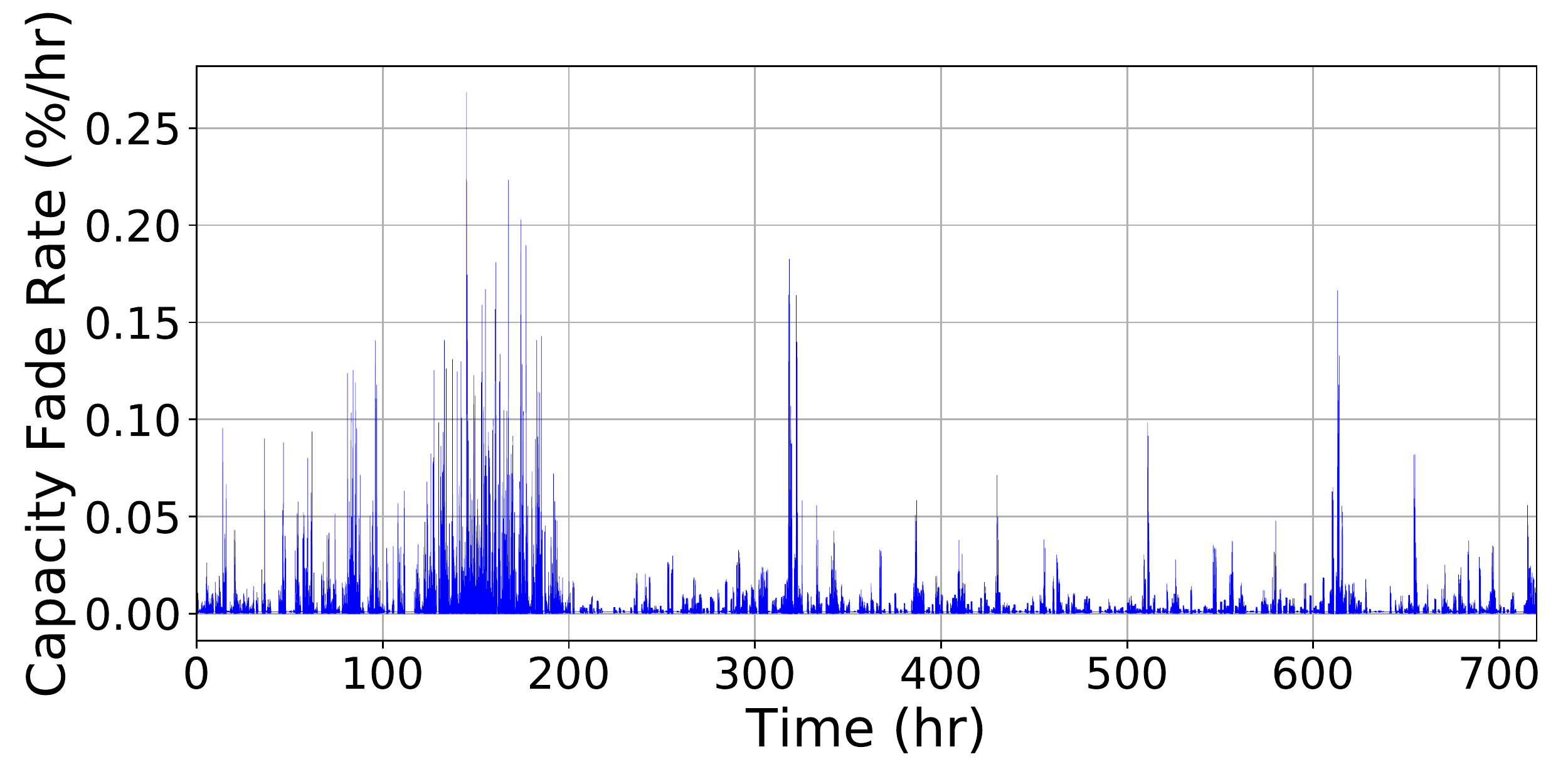}
\includegraphics[width=0.6\textwidth]{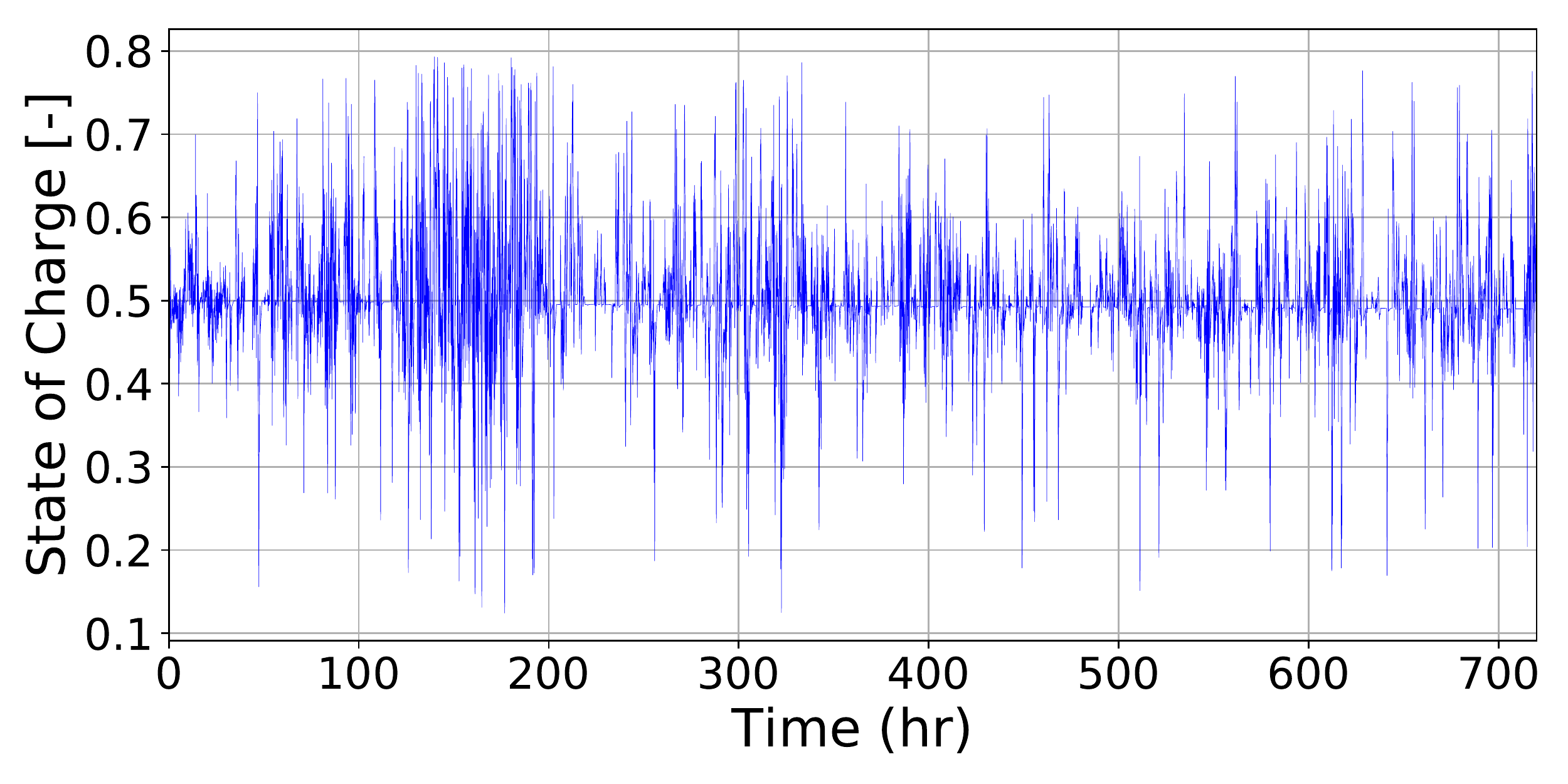}
\caption{Time profiles for battery internal state variables for the first month using HF-MPC.}
\label{fig:MPC_states_first_month}
\end{center}

\end{figure}

Figure \ref{fig:MPC_decisions_first_month} shows time profiles for commitment variables including FR capacity committed and the amount of power purchased in the day-ahead-market for the first month. The corresponding price and FR signal data are shown in Figure \ref{fig:data}. Overall, we can see that HF-MPC allocates the FR band more conservatively. Specifically, only 20\% of the FR band committed is larger or equal to 3 MW. This MPC formulation only allocates the FR band aggressively when the FR price is favorable (e.g., between hours 150-200). 
Figure \ref{fig:MPC_states_first_month} shows the profiles of selected state variables (including capacity fade and state of charge). When the committed FR band is large (e.g., between hours 150-200) capacity fade increases at a high rate. This illustrates how {\em aggressive FR market participation affects battery internal states}. 

Table \ref{tb:comparison_year} summarizes the performance of the different strategies. We observe that LF-MPC improves the revenue of the fixed band policy by \$62,000. Most of the improvement is due to the reduction in cost (less power is purchased). Remarkably, {\em HF-MPC increases the lifetime of the battery by 143\%} (compared with LF-MPC), increases the cumulative FR band by 10\%, and {\em increases profit by 35\%} (\$107,000).  From this table we also observe that increasing the prediction horizon of LF-MPC decreases profit (this is inconsistent with typical MPC formulations). We attribute this inconsistent behavior to long-term battery degradation effects that the LF-MPC controller does not  account for. 

\begin{table*}[!h]
\caption{ \bf Comparison of fixed band, low-fidelity MPC, and high-fidelity MPC strategies over two years of operation.}
\label{tb:comparison_year}
\centering
\footnotesize
\begin{tabular}{cccccccc}
\hline	
      & Horizon & Life Time  &Revenue & Cost   & Profit  &Cumulative  & Purchased Power \\
        & (h)       & (days)      & ($\times$\$1,000)        &   ($\times$\$1,000)    &   ($\times$\$1,000)         &FR band (MW)     &  (MWh) \\
\hline						   								  
\begin{tabular}{@{}c@{}}Fixed FR band \\  (3MW) \end{tabular}   &1   &   156  &  358    & 118     &240    & 8972  & 885 \\
\hline		
\begin{tabular}{@{}c@{}}Fixed FR band \\  (10MW) \end{tabular}   &1   &  87   &   338   &    97  & 241    & 7803&  702   \\
\hline
\begin{tabular}{@{}c@{}}Low-Fidelity  \\ MPC  \end{tabular}   &1   &  128   & 361     & 58   & 303    &   8273&  447  \\
\hline
\begin{tabular}{@{}c@{}}Low-Fidelity \\ MPC  \end{tabular}   &24   & 113    &352      &55      &297   &   7900&   423   \\
\hline
\begin{tabular}{@{}c@{}}High-Fidelity  \\ MPC  \end{tabular}   &1   &  312   & 474     & 64     &410    & 9061      & 490   \\
\hline						   							 								 	
\end{tabular}
\end{table*}

\subsection{Modification of Low-Fidelity Model}

The cumulative FR band is the total amount of FR band committed throughout the battery lifetime and can be viewed as an {\em effective lifetime}. Our closed-loop simulations indicated that, compared with LF-MPC with a one-hour prediction horizon, HF-MPC improves the cumulative FR band by 10\% and improves profit by 35\%.  We hypothesize that these improvements are mainly due to the long-term capacity fade effect (captured in the penalty term in the objective function). This term forces HF-MPC to allocate FR capacity {\em only when short-term market conditions are favorable relative to the long-term capacity effect}. In other words, the penalty ${\pi}^{C^f}$ in HF-MPC represents the long-term battery value.  The heuristic and LF-MPC strategies do not capture this long-term economic behavior. To verify our hypothesis, we modified the LF-MPC formulation; here, we assumed that capacity fade is a function of the FR band committed and we thus introduce a dynamic equation of the form $ C^f_{t+1,1} - C^f_{t,1} = \lambda\cdot  F^b_{t}$, where $\lambda$ is the percentage of capacity fade per MW of FR band. This function is introduced in the objective function \eqref{Eq:objective} to capture long-term capacity fade explicitly. Although this strategy is simple and neglects the dynamics of the FR signals and of the SEI, we will see that the performance of LF-MPC drastically improves. This is because the effective charging/discharging rate is moderate. Although we set the nominal maximum charging/discharging rate to be 10 MW, the effective charging/discharging rate $P$ is bounded by the maximum FR band times the FR signal $\alpha_{t,s}$. Our closed-loop simulations show that the charging/discharging rate remains below 3 MW for 99\% of the time when HF-MPC is employed,  and 92\% of the time when the maximum band policy (fixed policy at 10 MW) is used. Based on the nominal cumulative FR band value of $8,200$ MW, we estimate a value of $\lambda=\frac{20}{8200}=0.0024$ MWh/MW. 

Table \ref{tb:variant_simplified} summarizes the performance of LF-MPC  using a penalty term on capacity fade. We can see that this approach significantly improves the lifetime (by 119\%) and improves profit (by 26\%) over the original LF-MPC policy. The cumulative FR band is decreased by 28\%, which means that the battery allocates FR capacity more conservatively. For the modified LF-MPC formulation, increasing the prediction horizon from one to 24 hours improves the operational revenue (which is consistent with behavior of typical MPC formulations).  This consistency reinforces our observation that long-term economic effects of degradation indeed drive the policy of the controller and that such effects can be captured using a simple model. Although the performance of the modified LF-MPC formulation is still inferior to that of HF-MPC (in terms of profit), the performance gap is significantly reduced. These results highlight how one can use insights from detailed physical models to create improved MPC formulations of low computational complexity. In particular, the improved LF-MPC formulation is still a linear program that can be solved over a horizon of 24 hours and at high time resolutions (while the HF-MPC counterpart can only be solved for a 1-hour horizon). 

\begin{table*}[!h]
\caption{ \bf Performance of a modification of low-fidelity MPC.}
\label{tb:variant_simplified}
\centering
\begin{tabular}{cccccccc}
\hline	
       Horizon & Lifetime  &Revenue & Cost   & Profit  &Cumulative  & Purchased Power \\
         (h)       & (days)      & (\$1,000)        &   (\$1,000)    &   (\$1,000)         &FR band (MW)     &  (MWh) \\
         \hline	
 1   &  307   & 411     & 42   & 369    &5472    & 287   \\
\hline
  24   & 259    &434     &41      &393   &5708  & 288    \\
\hline			   							 								 	
\end{tabular}
\end{table*}

\subsection{Impact of Flexible Load and Capacity Fade Value}

We have also used HF-MPC to explore effects of flexibility gained by adjusting the load {\em every 2 seconds} (instead of every hour). This provides more degrees of freedom to the control formulation. Table \ref{tb:flex_load} shows that using a flexible load makes HF-MPC more aggressive in allocating FR capacities, which translates into a shorter lifetime and a higher cumulative FR band. The cost doubles due to more power ordered from the DAM, which is justified because the increase in revenue exceeds the rise in cost. Overall, however, a flexible load can improve operational profit by \$36,000 (8.8\%).   This illustrates how strategic manipulation of loads can help balance battery lifetime and overall profit. 

\begin{table*}[!h]
\caption{ \bf Effect of flexible load on performance of high-fidelity MPC.}
\label{tb:flex_load}
\centering
\begin{tabular}{cccccccc}
\hline	
      & Lifetime  &Revenue & Cost   & Profit  &Cumulative  & Purchased Power \\
             & (days)      & ($\times$\$1,000)         & ($\times$\$1,000)     & ($\times$\$1,000)          &FR band (MW)     &  (MWh) \\
\hline						   								  
\begin{tabular}{@{}c@{}} Constant  \\ Load  \end{tabular}   &  312   & 474     & 64     &410    & 9061      & 490   \\
\hline
\begin{tabular}{@{}c@{}}Flexible  \\ Load  \end{tabular}     &  271   & 572     & 126     &446    & 11444      & 1020   \\
\hline						   							 								 	
\end{tabular}
\end{table*}
\begin{table*}[!h]
\caption{\bf Effect of capacity fade value on high-fidelity MPC.}
\label{tb:capacity_fade_penalty}
\centering
\begin{tabular}{cccccccc}
\hline	
  ${\pi}^{C^f}$    & Lifetime  &Revenue & Cost   & Profit  &Cumulative  & Purchased Power \\
    ($\times$\$1,000/\%)         & (days)       & ($\times$\$1,000)         & ($\times$\$1,000)      & ($\times$\$1,000)         &FR band (MW)     &  (MWh) \\
\hline						   								  
 12  &  312   & 474     & 64     &410    & 9061      & 490   \\
\hline
20.5   &  410   & 510     & 58     &452    & 8771      &440    \\
\hline	
22.6     &426   &563 &59 & 504      &8502 &422 \\
\hline						   							 								 	
\end{tabular}
\end{table*}

The parameter ${\pi}^{C^f}$ represents the {\em long-term valuation} of battery capacity. As expected, the choice of this parameter is essential as it trades-off short-term and long-term economics. If the value of ${\pi}^{C^f}$ is too small, the effect of capacity fade is neglected and the HF-MPC controller will be more aggressive in allocating FR capacity. On the other hand, if ${\pi}^{C^f}$ is too large, the controller will become conservative in participating in the market. In the previous results, we set ${\pi}^{C^f}$ to 12,000\$/\%. This value was estimated as the profit of the fixed band policy (\$240,000) divided by 20\% (the remaining capacity at the end of life). Using this value, the HF-MPC controller achieved a profit of \$410,000. To analyze the effect of  ${\pi}^{C^f}$, we increased its value to 20,500\$/\% (obtained by dividing $\$410,000$ by 20\%). Table \ref{tb:capacity_fade_penalty} summarizes the results. We see that a larger value of ${\pi}^{C^f}$ makes HF-MPC more conservative in allocating FR capacities. Specifically, a longer lifetime and a smaller cumulative FR capacity are obtained. For this case, the controller increases profit by \$42,000 (10\%).  A more rigorous determination of the long-term value of capacity fade is an interesting topic of future work. 

\section{Conclusions}

We have presented a multiscale MPC framework to manage short-term economic value obtained from market transactions (energy and frequency regulation) and long-term economic value due to battery degradation. Insights gained from detailed closed-loop simulations provided insights to construct a low-complexity MPC formulation that can capture multiscale effects. \mynote{We believe that our proof-of-concept results can be of industrial relevance, as vendors are seeking to use batteries to participate in fast-changing electricity markets while maintaining asset integrity \cite{kumar2018stochastic,kumar2019benchmarking}}. As part of future work, we will seek to incorporate \mynote{more detailed battery models such as the Doyle-Fuller-Newman model, stochastic MPC formulations to capture market uncertainty,} and we will seek to accelerate simulations using parallel computers and machine learning techniques. 

\section*{Acknowledgments}

VZ acknowledges funding from the National Science Foundation under award NSF-EECS-1609183 and from the U.S. Department of Energy under grant DE- SC0014114. Work by VS at the University of Texas at Austin was supported by U.S. DOE award DEAC05-76RL01830 through PNNL subcontract 475525.  VZ declares a financial interest in Johnson Controls International, a for-profit company that develops control technologies. 

\appendix

\section*{List of Symbols}

\begin{footnotesize}
\noindent Battery
\setlist[description]{font=\normalfont}
\begin{description}[align=right,labelwidth=2cm, noitemsep,nolistsep]
\item [$D_j$] solid phase diffusion coefficient of lithium in the particles of electrode $j$, $m^2/s$
\item [$E$] remaining energy of battery, MWh
\item [$E_{max}$] battery capacity, WWh
\item [$F$] Faraday constant, 96,487 $C/mol$
\item [$c_e$] concentration of electrolyte in solution phase, $mol/m^3$
\item [$c_j$] solid phase concentration of lithium in the electrode $j$, $mol/m^3$
\item [$c^{avg}_j$] average concentration within the particle in the electrode $j$, $mol/m^3$
\item [$c^s_j$] surface concentration of lithium in the electrode $j$, $mol/m^3$
\item [$C_r$] rate of capacity fade, $1/s$ 
\item [$C_f$] capacity fade
\item [$i_{o,sd}$] the exchange current density for the side reaction $A/m^2$
\item [$I_{app}$] applied current passing through the cell, $A$
\item [$j$] negative $n$ and positive $p$ electrodes
\item [$J_j$] local reaction current density referred to electroactive surface area of electrode $j$, $A/m^2$
\item [$k_j$] rate constant of electrochemical reaction, $A/m^2/$, $m^{2.5}/(mol^{0.5}  s)$
\item [$M_{sd}$] molecular weight of SEI, $kg/m^3$
\item [$r$]  radial coordinate in the spherical particle, $m$
\item [$P$]  power supplied to the battery, MW
\item [$R_j$]  particle radius for electrode $j$, $m$
\item [$R_f$]  resistance of the SEI film, $\Omega \cdot m^2$
\item [$R_{SEI}$]  initial resistance of the SEI film, $\Omega \cdot m^2$
\item [$sd$] side reaction
\item [$S_j$]  total electroactive surface area of electrode $j$, $m^2$
\item [$t$]  time, $s$
\item [$T$] temperature, $K$
\item [$U_j$] local equilibrium potential of electrode j, $V$
\item [$U_{ref}$] constant equilibrium potential of the side reaction, $V$
\item [$V$]  cell voltage, $V$  
\item [$\phi_j$] solid phase potential of electrode $j$, $V$ 
\item [$\kappa_{sd}$] SEI ion conductivity, $S/m$ 
\item [$\eta_j$]  local overpotential of electrode $j$, $V$ 
\item [$\rho_{sd}$] density of SEI, $kg/m^3$
\item [$\delta_f$] SEI film thickness, $m$
\end{description}

\noindent Electricity Market
\setlist[description]{font=\normalfont}
\begin{description}[align=right,labelwidth=2cm, noitemsep,nolistsep]
\item [$F_k$] FR capacity provided at $k$-th hour, MW
\item [$L_k$]  committed power to load at $k$-th hour,  MW
\item [$N$] length of the prediction horizon
\item [$O_k$] power purchased from the day-ahead-market at $k$-th hour, MW
\item [$P_{k,s}$] net battery charge/discharge rate at $k$-th hour and $s$-th step, MW 
\item [$\overline{P}$] maximum charging rate, MW
\item [$\underline{P}$] maximum discharging rate, MW
\item [$S$] number of time steps per hour
\item [$x_{k,s}$] state variables of battery at $k$-th hour and $s$-th step 
\item [$\alpha_{k,s}$]  fraction of FR capacity requested by the ISO at $k$-th hour and $s$-th step 
 \item [$\eta_l$] minimum percentage of energy at the end of the prediction horizon
\item [$\eta_u$] maximum percentage of energy at the end of the prediction horizon
\item [${\pi}^{C^f}$] penalty parameter for capacity fade \$/\%
\item [$\pi^e_k$] electricity price at $k$th hour, \$/MWh
\item [$\pi^f_k$] FR capacity price at $k$th hour, \$/MW
\item [$\tau_l$] minimum percentage of energy within battery
\item [$\tau_u$] maximum percentage of energy within battery
\item [$\Phi$] total amount of profit earned before the end of life

\end{description}

\end{footnotesize}

\bibliographystyle{elsarticle-num}
\bibliography{first_principles.bib}

\end{document}